
\input amstex
\documentstyle{amsppt}

 \magnification=1200
 \NoBlackBoxes

 \TagsOnRight
 \topmatter
 \title
 Localisation of the Donaldson's invariants ...
 \endtitle
 \title Localisation of the Donaldson's invariants along Seiberg-Witten
classes.
 \endtitle
 \author Victor Pidstrigach and Andrei Tyurin
 \endauthor
 \author Victor Pidstrigach and Andrei Tyurin
 \endauthor
 \address  Steklov Mathematical Institute, ul. Vavilova 42, Moscow, 117966,
GSP-1, Russia
 \newline
and
\newline
Mathematical Institute, 24-29 St.Giles, Oxford OX1 3LB, UK
\newline
 e-mail: pidstrig\@class.mian.su, pidstrig\@maths.ox.ac.uk ,
tyurin\@tyurin.mian.su
 \newline
 \endaddress

 \abstract
 This article is a first step in establishing a link between the
Donaldson polynomials and Seiberg-Witten invariants of a smooth
4-manifold.

 \endabstract
 \endtopmatter

 \define\proj{\Bbb P} 

 \define\sA{{\Cal A}} 
 \define\sB{{\Cal B}} 
 \define\sC{{\Cal C}} 
 \define\sD{{\Cal D}} 
 \define\sG{{\Cal G}} 
 \define\sI{{\Cal I}} 
 \define\sM{{\Cal M}} 
 \define\sN{{\Cal N}} 
 \define\sO{{\Cal O}} 
 \define\sP{{\Cal P}} 
 
 
 \define\sT{{\Cal T}} 
 

 \document

\head
\S 0. Introduction.
\endhead

 As the result of analysis of the structure of the Donaldson's invariant of a
smooth simply
connected 4 manifold $X$ Kronheimer and Mrowka [K-M1] (see also [F-S]) proved
the  existence of
$k_X$ classes
 $$
K_1, \dots, K_{k_X} \in H_2(X, \Bbb Z)
 $$
and $k_X$ polynomials
$$
f_1, \dots, f_{k_X} \in \Bbb Q [z]
$$
such that the Donaldson's invariant of $X$ is defined uniquely by these pairs.

These classes are called the basic classes of $X$ and
the set $\{ K_i \},\text{for}~~~ i= 1, \dots, k_X, $ in
$H^2(X, \Bbb Z)$
is a diffeomorphism invariant.
So the diffeomorphism group $Diff X$ admits a representation to
the symmetric group on $k_X$ letters. A little bit later Seiberg and Witten
introduced
another set of classes $\{ SW \} $ and polynomials and Witten predicted the
shape of the
Donaldson's invariant in terms of these data. S-W classes are labelling the
moduli spaces
$\{ \sM_{SW} \}$ of the solutions of the Seiberg-Witten system of equations
(see \S1) and
 the dimensions of these moduli spaces are the degrees of the polynomials.
Using other type
 of system of equations we will prove that the classes of Kronheimer and Mrowka
and
Seiberg-Witten coincide. This paper is the first step of this proof and
contains those
facts from  gauge theory which we need for the proof.

At various stages of this project useful were conversations with Ed Witten,
S.Donaldson, S.Bauer.
We are grateful to S.Donaldson who invited us to the conference at the Newton
Institute in December 1994 where we described these results and to
I. Hambleton for his kind invitation to MPI, Bonn.
First author aknowldges with gratitude financial support of AMS (fSU grant)
and VW-stiftung during the work over this project as well as the kind
invitation
of S.Bauer and the hospitality of Bielefeld University.

\head
\S 1. Configuration spaces and its cohomology.
\endhead

 Let $c \in H^2(X)$ be a $Spin^c$ structure on X, i.e. a pair of
 complex hermitian rank 2 bundles $W^{\pm} _c, c_1(W^{ \pm} _c) =c$ (we
shall omit subscript $c$ if there is no confusion). Then the tangent bundle
$$
TX = Hom_h (W^-, W^+)
\tag1.1
$$
is the bundle of homomorphisms preserving the Hermitian structures.

 Fix a pair of connections $ \nabla^{\pm} $
 on $W^{\pm}$ such that their tensor product is Levi Civita connection
on $TX$
(which is the same as fixing the determinant $ \nabla_{det} =det \nabla^{\pm}
$).
This defines a Dirac operator:
$$
\sD^{g, \nabla_{det}} : \Gamma (W^+) \to \Gamma (W^-)
\tag 1.2
$$
as a composition of the connection $ \nabla^+: \Gamma (W^+) \to \Gamma (W^+
\otimes T^*X)$ and the convolution along (1.1). We shall omit the superscripts
if this will not make a confusion.

\subheading{U(1) - monopoles}

 We shall remind that a configuration space $\sC_{det W^+}$
 for $U(1)$-Seiberg-Witten moduli space is
 $$
 \sA_{det W^+} \times \Gamma (W^+) / \sG_{det W^+} =
 \sN \times \Gamma (W^+)/S^1
 \tag 1.3
 $$
 where $\sN = ker (d^* : \Omega ^1 \to \Omega^0)$ is a slice
 of an action of $\sG_{det W^+}$.
Seiberg-Witten moduli spaces $\sM_{SW}(0,c)$ is defined as a space of solutions
of the
system

$$
\sD^{ \nabla_{det}} \phi = 0
$$
$$
F_{\nabla_{det}}^+ = -(\phi \otimes \overline {\phi})_0,
\tag 1.4
$$
 modulo action of the gauge group,
where $( \nabla_{det}, \phi) \in  \sA_{det W^+} \times \Gamma (W^+)$.
Over the subspace $\sC^* _{detW^+} = \sN \times \Bbb P( \Gamma (W^+))$
 of non-reducible points of  $\sC_E$ there is a complex linear
 bundle $\sO_{\sC^*}(1) = \sO _{ \Bbb P( \Gamma (W^+))} (1)$
with the first Chern class $t=c_1( \sO_{\sC^*}(1) )$ generating
cohomology ring $H^*(\sC^*, \Bbb Q)$ in the case $X$ is simply
connected.

The system of equations (1.4) can be gauge invariant perturbed as
$$
\sD^{ \nabla_{det}} \phi = 0
$$
$$
(\omega + F_{\nabla_{det}})^+ = -(\phi \otimes \overline {\phi})_0,
\tag 1.$4_{\omega}$
$$
where $\omega$ is any 2-form.

Seiberg-Witten moduli spaces $\sM_{SW}(\omega,c)$ is defined as a space of
solutions of the
perturbed system (1,$4_{\omega}$)  modulo action of the gauge group.
For the generic form $\omega$ the moduli space is cut out transversally and
therefore smooth.

We shall also consider another perturbation - perturbation of the metric $g$.
The similar
"generic position" statement for this space of parametra is provided by the
following

\proclaim{Lemma 1.1}Fix some 2-form $\omega$ and consider the map
$$
\sM etr \times \Omega^1 \times \sA_{det W^+} \times \Gamma (W^+) @>f>>  \Gamma
(W^-)\times \Omega^2_+
\tag1.5
$$
given by sending
$$
(g, \xi, \nabla_{det}, \phi) \to (\sD^{g, \nabla_{det}} \phi,(\omega +
F_{\nabla_{det}})^+ + (\phi \otimes \overline {\phi})_0.
\tag1.6
$$
This map is a submersion.

\endproclaim

\demo{Proof}We shall remark  first that since principal symbol of Dirac
operator, as well as of Levi-Civita connection, doesn't depend on metric,
dependence is given by the formula:
$$
\sD^{g+r, \nabla_{det}} = \sD^{g, \nabla_{det}} + \delta_r
\tag1.7
$$
where $\delta_r \in \Omega^1$. Thus
$$
\delta f (r, \alpha, \delta \nabla_{det}, \delta r) = (\delta_1 f, \delta_2 f)
$$
where
$$
\delta_1 f = (\alpha +  \delta\nabla_{det} + \delta_r) \otimes \phi +  \sD^{g,
\nabla_{det}}(\delta\phi),
$$
$$
\delta_2f = r^*(F_{\nabla_{det}}^+ + \omega) +
d_{\nabla_{det}}^+(\delta\nabla_{det}) + (\delta\phi \otimes \overline
{\phi})_0 + (\phi \otimes \overline {\delta\phi})_0.
\tag1.8
$$
It is the standard observation that the restriction
$$
\delta_1f_{ \vert r = \delta\nabla_{det} = 0} = (\alpha + \delta_r) \otimes
\phi +  \sD^{g, \nabla_{det}}(\delta\phi)
$$
is an epimorphism. Therefore it is enough to prove that the restriction
$$
\delta_2f \vert_{\alpha = \delta\phi = 0}
$$
is onto. This follows from the fact that for any non zero 2-form
$(F_{\nabla_{det}}^+ + \omega)$ and for any harmonic self dual 2-form $\chi$
one can find an infinitesimal deformation of the metric $r$ such that the
pairing
$$
\langle r^*(F_{\nabla_{det}}^+ + \omega), \chi \rangle \neq 0.
$$
Indeed, by Cor.4.2.23 of [D-K] such $\chi$ vanishes on an open set if only
$\chi = 0$ everywhere.

\enddemo

There is another gauge invariant perturbation of the equation. Let
$$
\sT_{\frac 1 2} \subset \Gamma (End(\Lambda^0 \oplus \Lambda^2 _+))
\tag1.9
$$
be a subspace
of endomorphisms $t$ of the bundle $\Lambda^0 \oplus \Lambda^2 _+$
subject to the condition $\vert t \vert \le \frac 1 2 $.
Any such $t$ defines equation as follows
$$
\sD^{ \nabla_{det}} \phi = 0
$$
$$
(\omega + F_{\nabla_{det}})^+ = ((t-1)(\phi \otimes \overline {\phi}))_0,
\tag 1.$4_{\omega,t}$
$$

At points where $ \vert \phi \vert $ has its maximal value using Weitzenb\"ock
formula one has
$$
0 \le \Delta  \vert \phi \vert \le \frac {-s} 2  \vert \phi \vert ^2
+ \langle ((t-1)(\phi \otimes \overline {\phi}))_0 (\phi), \phi \rangle
=
$$
$$
 = \frac {-s} 2  \vert \phi \vert ^2 - \vert \phi \vert ^4
+ \langle (t(\phi \otimes \overline {\phi}))_0 (\phi), \phi \rangle
\le
\frac {\vert \phi \vert ^2} 2 ( -s - \vert \phi \vert ^2 + \vert t \vert
 \vert \phi \vert ^2 )
\tag1.10
$$
and as a consequence
$$
\vert \phi \vert  ^2 \le \frac {-s} { (1- \vert t \vert)} \le - 2 s
$$
provided our condition $\vert t \vert \le \frac 1 2$ is true
(compare [K-M2]).

Thus following well known properties of  $\sM_{SW}(\omega,c)$  for generic
metric $g$ (see for example [K-M2], [W]) are true also for the moduli
$\sM_{SW}(\omega, t,c)$ of solutions of 1.$4_{\omega,t}$:

1) there exists only finite set of classes $ c \in H^2(X)$ such that
$\sM_{SW}(\omega, t,c) \neq \emptyset$;

2) every moduli space  $\sM_{SW}(\omega, t,c)$ is compact;

3) if $X$ admits a metric with positive scalar curvature then
$\sM_{SW}(0, t,c) = \emptyset$ (in a relevant chamber in the case
$b_2 ^+ (X) = 1$).

So we can see that the Seiberg-Witten system of equations (1.4) and
(1.$4_{\omega, t}$) is very strong and the space of solutions is very rigid.

The strategy to get a non trivial solution is to consider a degenerated
solution
$(\nabla_{det}, 0) \in \sM_{SW}(\omega,c)$ for some $\omega$ and to use the
description of Kuranishi map for this singular point.

\subheading{Non abelian monopoles}

Now take a $PU(r)$-bundle $\xi$ on $X$ and fix its lift to a $U(r)$-bundle $E$
that is
fix $c_1(E) \in H^2(X)$. As a configuration space consider the set of triples
$$
(a_0, \nabla_{det}, \phi) \in \sA(P(E)) \times \sA(det(E \otimes W^+))
\times \Gamma (E \otimes W^+)
\tag1.11
$$
where $a_0$ is a traceless part of $U(r)$-connection, $\nabla_{det}$ is an
abelian connection on the line bundle $ det(E \otimes W^+)$ and $\phi$ is a
section of the vector bundle $E \otimes W^+$. The first two components of a
triple (1.11) define a coupl
ed Dirac operator
$$
\sD^{ \nabla_{det}}_{a_0}: \Gamma (E \otimes W^+) \to \Gamma (E \otimes W^-)
\tag1.12
$$
by the Leibnitz rule for  the action of the connection along $E$ and the
ordinary Dirac operator (1.2) along $W^+$. Now we can consider the non abelian
analogy of the system (1.4) (see [V-W]) :
$$
\sD^{ \nabla_{det}}_{a_0}(\phi) = 0
$$
$$
(F_{a_0} + \frac{1}{2} F_{\nabla_{det}} \otimes id_E)^+ = -(\phi \otimes
\overline {\phi})_0,
\tag 1.13
$$
and   the symmetry group of this
  system is given as the central extension
$$
1 \to \sG_{det(E\otimes W^+)} \to \sG \to \sG_{\xi} \to 1.
\tag1.14
$$
The moduli spaces $\sM_{SW}(0,2c_1 + rc, p_1)$ of Seiberg-Witten monopoles are
defined as a space of solutions of this system   modulo action of the gauge
group (1.14), where $p_1$ is the Pontriagin number of $E$.

In the same vein we can define the $\omega, t$ perturbed system considering
$$
(F_{a_0} + ( \frac{1}{2} F_{\nabla_{det}} + \omega) \otimes id_E)^+ =
((t-1)(\phi \otimes \overline {\phi}))_0
\tag1.1$3_{\omega}$
$$
as the second equation of the system (1.13) and get the moduli space
$$
\sM_{SW}(\omega, t,2c_1 + rc, p_1).
$$

In particular
$$
\sM_{asd} \subset \sM_{SW}( -\frac{1}{2} F_{\nabla_{det}}, t, 2c_1 + rc, p_1)
\tag1.15
$$
as the subspace of degenerated solutions of type  $(a_0,  \nabla_{det}, 0)$.

Non abelian Seiberg-Witten equations are very strong too (in spite of fact that
we lost a compactness). In particular it is easy to see that

1) in the symplectic (or K\"ahler) case we have nontrivial solution of (1.13)
only if
$$
(2c_1 + rc) \cdot [\omega] \leq 0
\tag1.16
$$
where $ [\omega]$ is the class of the symplectic (K\"ahler) form;

2) if $p_1$ is negative enough then generic degenerated solution $(a_0,
\nabla_{det}, 0)$ can be deformed to non trivial solution of the system
(1.1$3_{\omega}$) for any $\omega$.

3)  if $X$ admits a metric with positive scalar curvature then  $
\sM_{SW}(0,2c_1 + rc, p_1) = \emptyset$ always.

To get more plastic moduli space we construct the system of equations which
haven't any analogy      in the abelian case.

 For simplicity we will consider the case of $rk 2$ bundle only and  use the
following trick: if we fix some connection $b_0$ on $det W^+$ then the first
pair of components (1.11) $(a_0, \nabla_{det})$ is given uniquely by a
$U(2)$-connection $a$ on $E$.
So instead the space of triples (1.11) we can consider the space of pairs
$$
(a, \phi) \in  \sA(E)  \times \Gamma (E \otimes W^+)
\tag1.17
$$
as a configuration space and $\sG_{E}$ as a gauge group.

 Let $\omega \in \Omega ^2 (X) $ be a two form in the cohomology
 class of $- \pi i c_1(E)$.
 Denote by $\sA_{\omega} \subset \sA_E$ a subspace of the space of all $ U(2)$
 connections $a$ on the bundle $E$ subject to the condition
 $$
 tr F_a = 2\omega ,
 \tag 1.18
 $$
 $tr : EndE \to \Bbb C$ is a trace. Traceless part of the curvature
 $F_a - \frac 1 2 tr (F_a) $ will be denoted as $(F_a)_0 = F_{a_0}$.

 A gauge group $\sG_E$ of $E$ acts on both spaces $\sA_E$ and
 $\sA_{\omega}$ with orbit spaces, resp., $\sB_E$ and $\sB_{\xi}$.

 Let us consider a natural action of this group on the space
 $\sA_{\omega} \times \Gamma (W^+ \otimes E)$ with the orbit
 space denoted as $\sC_E$. So $\sC_E$ is the fibration
$$
\sC_E \to \sB_{\xi}
\tag1.19
$$
with a fibre $\Gamma (W^+ \otimes E)$.

 Let $(a, \phi) \in \sA_{\omega} \times \Gamma (W^+ \otimes E)$.
 Denote by $N_{(a, \phi)}$ a slice of the action at $(a, \phi)$:
 $$
 N_{(a, \phi)} = Ker ( d_a ^* \oplus m_{\phi} ^*)
 \tag1.20
 $$
 where
 $$
 d_a \oplus m_{\phi}: \Omega ^0 (ad(E)) \to \Omega ^1 (ad(E))
 \times \Gamma (W^+ \otimes E),
 \tag 1.21
 $$
 is a tangent map of the action, $m_{ \phi}$ is a multiplication of a
 vector spinor $ \phi$ by the endomorphism.

\subheading{Stabilisers of the gauge action}

 The stabiliser $St = St_{(a, \phi)}$ of a point $(a, \phi)$ is at most
 $U(2)$ - the largest possible stabiliser of the connection itself.
 There are three possibilities for $St_a \subset U(2)$ (provided by the
 fact that $St_a$ is a centraliser of some subgroup of $U(2)$) except
 the trivial subgroup: $S^1, S^1 \times S^1, U(2)$. The last one
 corresponds to the pair consisting of trivial connection
 and zero spinor. Its neighbourhood is modelled as $N_{(0, 0)} / U(2)$.

 Centre $S^1 \subset U(2)$ of the structure group
 is the stabiliser of a point of
 the type $(a, 0)$, where $a$ is non-reducible. Orbits of such a
 points form the subspace of the zero-section of the fibration (1.19)
 $\sB_{\xi} \subset \sC_E$ with the neighbourhood $U = \tilde{U}/S^1$
 where $\tilde{U}$ is modelled on the total space of the bundle (1.19):
 $$
 \tilde{U} = \Gamma(W^+ \otimes E)  \times _{\sG_E /S^1} \sA^*_{\omega}
 \to \sB_{\xi}^*,
 \tag 1.22
 $$
 $S^1 $ acts with a weight one representation on fibres.
 The determinant of the universal bundle $det \Bbb E$  can be described as
 $S^1$ -equivariant bundle over $X \times \tilde{U}$ associated to a
 tautological representation of $S^1$ in $det E$.

 If the stabiliser $S^1$ of $(a,\phi)$ is a different subgroup of $U(2)$,
 namely the one which maps onto a maximal torus of $PU(2)$,
 it means that connection $a= \lambda_1 \oplus \lambda_2 $
 is reducible
 due to some reduction of the bundle $E = L_1 \oplus L_2$,
 and with spinor of the form $ (\phi_1  \oplus \phi _2), \phi _1 =0,
 \phi_2 \ne 0$ with respect to the same reduction of the bundle.
 Subspace of all orbits of such a singularities is a subspace of
 $$
 \sA_{L_1} \times
 \sA_{L_2} \times \Gamma (W^+ \otimes L_2) / (\sG_{L_1} \times \sG_{L_2})
 \tag 1.23
 $$
 consisting of points subject to the condition
 $ F_{\lambda_1} + F_{ \lambda_2} = 2 \omega$ (which,  given $\lambda_2$,
 determines $\lambda_1$ up to a gauge equivalence )
 i.e. isomorphic to $ \sC_{L_2}$.
 The neighbourhood $U$ of this locus is modelled as follows:

 $$
 U = \sN \times (\Omega^1 (L_1 \otimes L_2 ^{-1} ) \times
 \Gamma (W^+ \otimes L_1) \times
 (\Gamma (W^+ \otimes L_2)-\{0\}))/(S^1 \times S^1).
 \tag 1.24
 $$
 This is a total space of the real cone fibration

 $$
 (\Omega^1 (L_1 \otimes L_2 ^{-1} ) \otimes
 \sO  _{\sC _{det (W^+ \otimes L_2)} ^*}(-1) \oplus
 \Gamma (W^+ \otimes L_1) \otimes
 \sO ^* _{\sC _{det (W^+ \otimes L_2)} ^*}) /S^1
 \tag 1.25
 $$
 with weight one action of $S^1$ on fibres of the vector bundle,
 over $\sC^* _{det (W^+ \otimes L_2)}$.
 The restriction of the universal bundle to $U$ is given as a $S^1 \times S^1$
-equivariant
 fibration over $S^1 \times S^1$-space

 $$
 X \times
 \sN \times (\Omega^1 (L_1 \otimes L_2 ^{-1} ) \times
 \Gamma (W^+ \otimes L_1) \times
 (\Gamma (W^+ \otimes L_2)-\{0\}))
 \tag 1.26
 $$
 given by natural action of  $S^1 \times S^1$ on $L_1 \oplus L_2$.
 This is equivalent to $S^1$-equivariant bundle $L_1 \oplus L_2 \otimes
 \sO _{\sC _{det (W^+ \otimes L_2)} ^*}(1) $
 over the total space of the bundle

 $$
 (\Omega^1 (L_1 \otimes L_2 ^{-1} )
 \otimes \sO _{\sC _{det (W^+ \otimes L_2)} ^*}(-1) \oplus
 \Gamma (W^+ \otimes L_1) \otimes \sO _{\sC _{det (W^+ \otimes L_2)} ^*}).
 \tag 1.27
 $$

\subheading{Cohomology ring of $\sC_E ^*$}

 Let $(\sA_{\omega} \times \Gamma (W^+ \otimes E))^*$ denotes the subspace
 of all points in  $\sA_{\omega} \times \Gamma (W^+ \otimes E)$ with
 trivial stabiliser. Denote by $\sC_E ^*$ corresponding orbit space.
 There is a universal bundle $\Bbb E$ over the space $X \times \sC_E ^*$
 $$
 \Bbb E = E \times _{\sG_E} (\sA_{\omega} \times \Gamma (W^+ \otimes E))^* .
 \tag 1.28
 $$

 The cohomology ring $H^*( \sC_E ^*, \Bbb Q)$ is
 generated by 2-dimensional classes
$$
\mu (\Sigma) = \frac 1 4 p_1( \Bbb E) / [\Sigma],
$$

$$
 t=c_1( \Bbb E) /[pt]
\tag 1.29
$$
and 4-dimensional class $ \nu = \mu (pt) = \frac 1 4
 p_1( \Bbb E) / [pt]$ as it follows from the fibration
 $$
 \Bbb R _+ \times \Bbb {CP} ^{ \infty} \to
 (\sA_{\omega} \times \Gamma (W^+ \otimes E))^* / \sG_E
 \to \sA _{ \xi} ^* / \sG_{ \xi}.
 \tag 1.30
 $$

 Space
 $$
 \sP_E = (\sA_{\omega} \times
 \Bbb S (\Gamma (W^+ \otimes E))^* / \sG_E),
 \tag 1.31
 $$
where $\Bbb S ( )$ is the unit sphere bundle, is a deformational retract of
$\sC _E ^*$ and of the blow-up
 $\hat \sC _E$ (cf. \S 3 below).
Therefore an inclusion
$\sP _E @>i>> \hat \sC _E$
together with the projection
 $\hat \sC _E @>i>> \sP _E$
induces isomorphism in cohomologies. Thus the restriction
of the $\mu, \nu$ classes to the locus of reducibles of the first type
(as well as on it's link in modified definition)
coincides with those for $\sP_E$.
Restriction of the $t$ class is the generator of the
fibre $ \Bbb {CP} ^{ \infty}$ of the trivial fibration (1.30).

For the locus of reductions of the second type description (1.26-1.27)
of the restriction of the universal bundle to the neighbourhood
of this locus gives restriction of our generators.
Link of this locus is a projective bundle

$$
\Bbb P _{L_1, L_2} =
\Bbb P (\Omega^1 (L_1 \otimes L_2 ^{-1} ) \otimes \sO _{\sC^* _{L_2}}(-1)
\oplus \Gamma (W^+ \otimes L_1) \otimes \sO)
\tag 1.32
$$
with the restriction of the universal bundle given as

$$
L_1 \otimes \sO
_{\Bbb P (\Omega^1 (L_1 \otimes L_2 ^{-1} ) \otimes \sO _{\sC^* _{L_2}}(-1)
\oplus \Gamma (W^+ \otimes L_1) \otimes \sO)} (1)
\oplus L_2 \otimes \sO _{\sC^*_{L_2}}(1).
\tag 1.33
$$
A small computation shows
$$
\mu (\Sigma) = \frac 1 2 \langle l_1 - l_2 , \Sigma \rangle
(t + c_1( \sO
_{\Bbb P (\Omega^1 (L_1 \otimes L_2 ^{-1} ) \otimes \sO _{\sC^* _{L_2}}(-1)
\oplus \Gamma (W^+ \otimes L_1) \otimes \sO)} (1)))
\tag 1.34
$$
$$
\nu =
(t + c_1( \sO
_{\Bbb P (\Omega^1 (L_1 \otimes L_2 ^{-1} ) \otimes \sO _{\sC^* _{L_2}}(-1)
\oplus \Gamma (W^+ \otimes L_1) \otimes \sO)} (1)))^2.
\tag 1.35
$$

 \head
 \S 2.The equation.
 \endhead

 Let  $a \in \sA_{E}, \phi \in \Gamma (E \otimes W^+)$ and
$\nabla^{\pm} \in \sA_{W^{\pm}}$. On the tensor product of bundles
$ Hom(W^-, W^+) \otimes ad E$  connections $a$ and $\nabla^{\pm}$ define a
connection
which is
the tensor product of  the Levi - Civita connection and the traceless part
$a_0$ of $a$.
Dirac operator (1.12) coupled with a connection $a$ is defined as a map:
$$
\sD_{a_0} ^{ det( a \otimes \nabla^{\pm})} :
\Gamma (W^+ \otimes E) \to \Gamma (W^- \otimes E).
\tag2.1
$$
  Here the dependence on  traceless part $a_0$ and
$U(1)$-connection
$$\nabla_{det}=
det( a \otimes \nabla)
$$
on $det ( E\otimes W^+)$
will be exploited in
different way and this is reflected in notations.

Fix any 2-form $\omega$ in the cohomology class of $- 4 \pi i f = - 4 \pi i
(c_1+c)$ and
$t \in \sT_{ \frac 1 4 }$.
 Then our system of equations is:
$$
 F_{\nabla_{det}} = \omega,
$$
 $$
\sD_{a_0}^{\nabla_{det}} \phi = 0,
$$
 $$
F^+_{a_0} =((t -1)(\phi \otimes \overline { \phi}))_{00}
\tag2.2
$$
  where  $((1-t)(\phi \otimes \overline { \phi}))_{00}$ is a ``double''
traceless
  component of
$$ ((1-t)(\phi \otimes \overline { \phi})) \in
  sl (E) \otimes \Lambda ^2 _+ = su (E) \otimes su(W^+) \subset u (E)
  \otimes u (W^+).
$$

  So the symmetry group of this
  system is given as the central extension (1.14)
which we can identify as $ \sG_E$ fixing any connection $b_0$ on $det W^{\pm}$.
Under this identification the coupled Dirac operator (2.1) depends on $a$ and
the metric $g$ on $X$. So
$$
\sD_{a_0}^{\nabla_{det}} = \sD_a^{g, c + c_1(E)}.
\tag2.3
$$

\subheading{Definition 2.1} We shall denote moduli space of orbits of
  solutions of (2.2), modulo action of the gauge group as
  $\sM _{B} ^{g, \omega,t}(c+c_1(E), p_1(E))$.

 The infinitesimal properties of this system such as the deformation complex,
normal cones of singularities and so on are related closely to the properties
of the following system of equations
$$
 F_{\nabla_{det}} = \omega,
$$
 $$
\sD_{a_0}^{\nabla_{det}} \phi = 0,
$$
 $$
F^+_{a_0} = 0
\tag2.$2_0$
$$
which was investigated in the series of papers [P-T], [T1-T5] and [P].

The moduli space of orbits of
  solutions of (2.$2_0$), modulo action of the gauge group was denoted as
  $\Cal {MP}(X, p_1(E), c+c_1(E))$ (the moduli space of jumping pairs, see Def.
of \S1 of [P] or (1.1.26) of [P-T]). Sending a pair $(a, \phi)$ to
$a_0$ we get the map
$$
\pi : \Cal {MP}(X, p_1(E), c+c_1(E)) \to \sM_{asd}
\tag2.4
$$
to the moduli space of $SO(3)$-instantons and the image of this map
$$
\sM _1 ^{g, \omega}(c+c_1(E), p_1(E)) = \pi (\Cal {MP}(X, p_1(E), c+c_1(E)))
\subset  \sM_{asd}
\tag2.5
$$
is called the moduli space of jumping instantons.

On the set of solutions of the system (2.2)
  the gauge group action is not free. Possible stabiliser groups
  are $S^1, S^1 \times S^1, U(2)$. The last one occurs as the stabilisers
  group of the solutions with trivial connection (on the trivial bundle)
  and zero spinor only.

  $S^1$ may be a stabiliser of two types of solutions:

\subheading{Singularities of the first type} Solutions with $\phi = 0$.
This implies $(F^+ _a)_0 = F_{a_0} ^+= 0$
  i.e. (uniquely) corresponding $PU(2)$ connection $ (a)_0$ is
  antiselfdual. Stabiliser is the centre of $U(2)$. We shall call
these solutions reducible of the first type. The subspace of reducible
solutions of the first type is diffeomorphic to the moduli of $asd$
connection for the adjoint bundle $\xi$.

Formally the set of singularities of this type is all $\sM_{asd}$
but
let
$$
 \sM _{B} ^{g, \omega,t}(c+c_1(E), p_1(E))_0 \subset  \sM _{B} ^{g,
\omega}(c+c_1(E), p_1(E))
\tag2.6
$$
be the subset of solutions with $\phi \neq 0$. Then in $\sM_{asd}$ the subset
$$
Sing_1 = lim  \sM _{B} ^{g, \omega,t}(c+c_1(E), p_1(E))_0
\tag2.7
$$
is defined correctly as the subset of limits of solutions with non trivial
twisted spinor field $\phi$.

\proclaim{Proposition 2.1} For generic metric $g$ and generic form $\omega$
$$
Sing_1 = \sM _1 ^{g, \omega}(c+c_1(E), p_1(E))
\tag2.8
$$
\endproclaim

This Proposition is following immediately from the Transversality Theorem (see
the next section) and the description of the Kuranishi family at a pair $(a,
0)$.

\subheading{Remark}There exists the elementary proof of this statement using
the iterative procedure providing by the third equation of (2.2).

\subheading{Singularities of the second type}Let $Sing_2$ be the set of
solutions with reducible connections, which we lift to $U(2)$-connection $a$,
then the Levi-Civita $PU(2)$-connection on $W^+$ is lifted to $U(2)$-connection
$\nabla^+$ and
$$
F_{\nabla^+} + F_{det a} = \frac 1 2 \omega .
\tag2.9
$$
Moreover,
  $a = \lambda_1 \oplus \lambda_2$,
  due to some reduction of the bundle $E = L_1 \oplus L_2$,
  and with spinor of the form $ \phi = 0 \oplus \phi _2$ with respect to
  the same reduction of the bundle.
  For a reducible connection the harmonic spinor always splits in this way
  and we require vanishing of the first component.
  Stabiliser here is of the form $S^1 \times 1 \subset S^1 \times S^1 \subset
  U(2)$. Let $l_i = c_1(L_i)$, $\alpha = l_1 - l_2$ and $f = c + c_1 (E)$.
Spinor $ \phi_2 \in \Gamma (W^+ _c \otimes L_2) =  \Gamma (W^+ _{c+2l_2}) $
is a harmonic one:
$$
\sD^{(2 \lambda_2 \otimes det \nabla^* = \nabla_{det})} \phi _2 = 0 .
$$
Now
$$
F_a  = \left( \matrix F_{\lambda_1} & 0 \\
                      0 & F_{ \lambda_2} \endmatrix \right)
  = \left( \matrix  \omega - F_{det \nabla^+} - F_{\lambda_2} & 0 \\
                      0 & F_{ \lambda_2} \endmatrix \right),
$$
$$
(F_a)_0  =\frac 1 2 \left(  \matrix  F_{\lambda_1} - F_{\lambda_2} & 0 \\
                      0 & F_{ \lambda_2} -  F_{\lambda_1} \endmatrix \right)
$$
and
$$
-(\phi \otimes \overline { \phi}) = \left( \matrix
 0 & 0 \\ 0 & -(\phi_2 \otimes \overline { \phi_2})
\endmatrix \right)
$$
$$
-(\phi \otimes \overline { \phi})_0 = \frac 1 2 \left( \matrix
(\phi_2 \otimes \overline { \phi_2})  & 0 \\
0 & -(\phi_2 \otimes \overline { \phi_2})
\endmatrix \right).
$$
Therefore second our equation
$(F^+_a)_0 = ((t-1)(\phi \otimes \overline { \phi}))_{00}$
provides

$$
 F^+ _{2 \lambda_2 + det \nabla^+} = ((t-1)(\phi_2 \otimes \overline {
\phi_2}))_0
- \omega^+.
\tag2.10
$$
So a reducible solution of the second type produces a solution
to $U(1)$ - Seiberg Witten equation:
$$
( 2 \lambda_2 + det \nabla^+ , \phi_2) \in
\sM_{SW} (\omega^+, t,f - \alpha).
\tag2.11
$$
It is easy to see that this is one to one correspondence and we proved

\proclaim{Proposition 2.2}
$$
Sing_2 = \bigcup_{\beta \vert (f-\beta)^2 = p_1} \sM_{SW} (\omega^+, t, \beta).
\tag2.12
$$
\endproclaim

\subheading{ Remark} The very important and very interesting case is when $f =
0$.
Here we can choose $\omega = 0$ and we have the standard Seiberg - Witten
equation (1.4). In this case the diffeomorphisms group $Diff_X$ acts on the
moduli space  $\sM _{B} ^{g, 0}(0, p_1(E))$ and all its singularities.

  And finally $S^1 \times S^1$ occurs as a stabiliser of a solution
  consisting of an $asd$ reducible connections and vanishing spinor.

\subheading{Linearization }  Linearization of the the system (2.2) is elliptic
and, therefore, Fredholm:

  $$
  \sD_a \delta \phi + \delta a * \phi = 0
  $$
  $$
  d^+ _a ( \delta a ) =
((t-1)( \delta \phi \otimes \overline{\phi} + \phi \otimes
  \overline{\delta \phi}))_{00}
\tag2.13
 $$
  where $ \delta a \in \Omega ^1 (u(E)) , \delta \phi \in \Gamma (E \otimes
  W^+)$ are components of tangent vector at the point $(a, \phi)$.

  This linear map together with the linearization of the action of
  gauge group gives following deformation complex, which is homotopic to a
  direct sum of a deformation complex for moduli space of $asd$-connections
  and Dirac operator:

  $$
  \Omega ^0 (u(E)) @> \delta_1>> \Omega ^1 (pu(E))
  \times \Gamma (W^+ \otimes E) @> \delta _2 >>
  \Omega ^2 _+ (pu(E))
  \times \Gamma (W^- \otimes E),
\tag2.14
 $$
  with first map given by $ \delta _1= d_a \oplus m_{\phi} $
  and second by a matrix
  $$
  \delta_2 ( \delta a, \delta \phi)=
  (d_a ^+ ( \delta a) +  ((t-1) (\phi \otimes \overline {\delta \phi} +
  \delta \phi \otimes \overline { \phi}))_{00},
  \sD_a (\delta \phi)   + (\delta a) \cdot \phi )
\tag2.15
 $$

 Index of the complex is given by following formula:

  $$
  ind = - \frac {3} 2 p_1 -3(1+b_2 ^+ ) + \frac 1 2 (f ^2 - \sigma_X) -1
 \tag2.16
 $$
  where $p_1 = p_1(E)$, $f = c + c_1(E)$, $ \sigma_X$ - signature of $X$.

  Standard technique of Kuranishi models gives a model
  of a neighbourhood of the point $(a, \phi)$ as space of
  $St_{(a, \phi)}$ orbits in a neighbourhood
  of zero in the preimage of zero of certain real-analytic map
  $ H^1 _{(a, \phi)} @> \Psi >> H^2 _{(a, \phi)} , $
  where $H_{(a, \phi)} ^i$ are homologies of the deformation complex.

  For a reduction of the first type one has all cross-terms
  in deformation complex vanishing and it turns to be a direct
  sum of the deformation complex for the moduli of $asd$-connections
  and Dirac operator with stabiliser $S^1$ acting only on the
  spaces of sections of twisted spinors $\Gamma (W ^{\pm} \otimes E)$
  with weight one.

  Let $(a, \phi) = ( \lambda_1 \oplus  \lambda_2, 0 \oplus \phi_2) $
  be a reduction of the second type. Its deformation complex
  is a sum of the deformation complexes of $U(1)$-monopole equations
  and following complex:

  $$
  \Omega ^0 (L_1 \otimes L_2 ^{-1}) @> \delta_1 >>
  \Omega ^1 (L_1 \otimes L_2 ^{-1})
  \times \Gamma (W^+ \otimes L_1) @> \delta _2 >>
  \Omega ^2 _+ (L_1 \otimes L_2 ^{-1})
  \times \Gamma (W^- \otimes L_1),
\tag2.17
 $$
  with first map given by $ \delta _1= d_{\lambda} \oplus m_{\phi_2} $
  and second by a matrix
  $$
  \delta_2 ( \alpha, \delta \phi_1)=
  (d_{ \lambda} ^+ ( \alpha)
+((t-1) (\phi_2 \otimes \overline {\delta \phi_1} +
  \delta \phi_1 \otimes \overline { \phi_2})_0,
  \sD_{\lambda_1} ( \delta \phi_1)
  +(\alpha) \cdot \phi_2) .
\tag2.18
$$

\head
\S 3. Transversality.
\endhead

Let us return to the begining of the previous section. If we fix the
connections $\nabla^{\pm}$ on the spinor bundles $W^{\pm}$ compatible with the
Levi - Civita connection then the Dirac operator (2.1) will depend on
$U(2)$-connection $a$ and as the conf
iguration space we can consider the space $\sA_{\omega} \times
\Gamma(W^+ \otimes E)$ with condition (1.18) instead the first equation of
(2.2).
We will do it in this section.

{}From Kuranishi description it follows that the point $(a, \phi)$
is smooth if the obstruction space vanishes: $H^2 _{(a, \phi)} = 0$.
Our moduli space depends also on different continuous parameters
i.e. metric, connection on $det W^+$ etc.
It is the task of this section to prove that for generic parameters
one has  $H^2 _{(a, \phi)} = 0$ for all solutions  $(a, \phi)$.

First we shall discuss transversality at reducible points.
In order to get it we shall change configuration space of our system.
Following illustrates reasons in the case of reducible solutions of the first
type.

Reductions of the first type, i.e. $(a, 0)$ with $asd$-connection
$a$  occurs for a generic metric only if
$v.dim (\sM_{asd}) = -2p_1(E) -3(1+b_2 ^+) \ge 0$ (cf [FU], Thm 3.1).
It also follows that for generic metric asd-component of the deformation
complex has
vanishing second cohomology $H_a^2 = 0$. Therefore the only possibility
for $H^2 _{(a, \phi)} \neq 0$ is when $coker \sD _a \neq 0$.
Generally it is the case since the Chern class of the minus index bundle
$c_{ind \sD _a +1} (-Ind \sD _a)$ does not vanish. This is a topological
obstruction for vanishing of the obstruction space $H^2 _{(a, 0)}$ for all
$asd$
connections $a$.

The Kuranishi description is given by a map
$$
\Psi : H^1 _a \otimes ker \sD_a \to coker \sD_a
\tag3.1
$$
with $\Psi (\alpha, \sigma )=
\alpha * \sigma + O(( \alpha, \sigma ) ^3 )$
has vanishing differential.

Thus we shall change equation, i.e.
configuration space to get rid of this topological
obstruction for $H^2 _{(a, 0)}$ to vanish everywhere.
Take instead of  the space
 $\sC _E \sA_{\omega} \times  \Gamma (W^+ \otimes E)$
its blow up in the locus $\sA_{\omega, red} \times \{0\}$ of pairs consisting
of the reducible connection and zero spinor:
$$
\widehat{\sC _E} = \widehat{ \sA_{\omega} \times  \Gamma (W^+ \otimes E)}.
\tag3.2
$$

Gauge group acts in a natural way on the blow up.
Lift the equation to the blow up in a natural way.
For example in the neighbourhood of reducibles of the first type for
$$
(a, \phi, \psi) \in \sA_{\omega} \times  \Gamma (W^+ \otimes E) \times
\proj (\Gamma (W^+ \otimes E))
$$
take the restriction of solutions of the system
 $$
\sD_a \psi = 0
$$
 $$
(F^+_a)_0 = -(\phi \otimes \overline { \phi})_{00}
$$
 to $Y \subset  \Gamma (W^+ \otimes E)
\times \proj (\Gamma (W^+ \otimes E)) $ where  $Y$ is a blow up of $\Gamma (W^+
\otimes E)$ in zero
or equivalently the total space of the linear bundle
$ \sO _{\proj (\Gamma (W^+ \otimes E))} (-1)$.

\subheading{Definition 3.1} Denote by
$\widehat{ \sM _{B} ^{g, \omega}}(c+c_1(E), p_1(E))$ moduli of solutions in the
space (3.2)
modulo the action of the gauge group.

There is an obvious projection
$$
\pi : \widehat{ \sM _{B} ^{g, \omega}}(c+c_1(E), p_1(E))
\to  \sM _{B} ^{g, \omega}(c+c_1(E), p_1(E))
\tag3.3
$$
(compare 2.4).

Now the space of reductions of the first type $\pi ^{-1} (\sM_{asd})$ can be
identified with the moduli space of pairs $\Cal {MP}(X, p_1(E), c+c_1(E))$
(2.4).
In [P-T] it is proved that this moduli space is smooth
in the sense that second homology of corresponding deformation
complex vanishes.

Consider reducible solutions of the second type $\pi^{-1}(\sM _{SW}(\beta))$.
Let
$$
(a = \lambda_1 \oplus \lambda_2, \phi = 0 \oplus \phi _2) \in \sM _{SW}(\beta)
$$
and
$$
(a, \phi_2, \langle \theta, \delta \phi_1 \rangle ) \in \pi^{-1}(\sM _{SW}(
\beta)),
\tag 3.4
$$
where
$$
\langle \theta, \delta \phi_1 \rangle \in
\Bbb P(ker \delta_1 ^*) \subset
\Bbb P( \Omega ^1 (L_1 \otimes L_2 ^{-1})
\oplus \Gamma (W^+ \otimes L_1) ).
$$

The curvature of $\lambda = \lambda _1 - \lambda _2$ has form
$F_{\lambda} = \sigma \otimes u$ for two-form $\sigma$ and $u \in ker d_{a_0}$.

Deformations of the solution of the modified system is described by  a sum
of the deformation
complex of Seiberg-Witten moduli space at the point $(\lambda _2, \phi_2)$
and following picture for the deformations in the normal direction (compare
2.17-2.18):
take the bundle

$$
\CD
 ( \Omega ^2 _+ (L_1 \otimes L_2 ^{-1})
  \oplus \Gamma (W^- \otimes L_1)) \otimes \sO_{\Bbb P} (1)
\\
@VVV \\
\Bbb P( ker \delta_1 ^* ),
\endCD
\tag3.5
$$
and take the section of this defined by the map $\delta_2$,
  with $\delta_i$ given by the formulas
  $$
\delta _1= d_{\lambda} \oplus m_{\phi_2}
$$
  $$
  \delta_2 ( \theta, \delta \phi_1)=
  (d_{ \lambda} ^+ ( \theta)
+((t-1) (\phi_2 \otimes \overline {\delta \phi_1} +
  \delta \phi_1 \otimes \overline { \phi_2}))_0 ,
  \sD_{\lambda_1} ( \delta \phi_1)
  +(\theta) \cdot \phi_2) .
$$
Now directions where deformation is not obstructed are given as zero set of
this section denoted by $s_{\delta_2}$.

Dependence on metric  given through that of $\delta_2$
and dependence on 1-form $ \eta $
 given by replacing $ \omega$ by $ \omega + d \eta $
presents $\delta _2$ as a section of
$$
\CD
 ( \Omega ^2 _+ (L_1 \otimes L_2 ^{-1})
  \oplus \Gamma (W^- \otimes L_1)) \otimes \sO_{\Bbb P} (1)
\\
@VVV \\
\Bbb P( \Omega ^1 (L_1 \otimes L_2 ^{-1})
  \oplus \Gamma (W^+ \otimes L_1) / \delta_1 (\Omega ^0 (L_1 \otimes L_2
^{-1})))
\times \sM etr \times \Omega ^1
\endCD
$$
with linearization at the point (3.4) given by the differential
$$
((\Omega ^1 (L_1 \otimes L_2 ^{-1})  \times
\Gamma (W^+ \otimes L_1)))  \times T \sM etr \times \Omega^1
\to
$$
$$
\to  \Omega ^0 (L_1 \otimes L_2 ^{-1})  \times
  \Omega ^2 _+ (L_1 \otimes L_2 ^{-1})
  \times \Gamma (W^- \otimes L_1),
\tag3.6
$$
$$
 (\theta ,  \delta \phi_1, r, \eta) \mapsto ( \delta_1 ^* ( \theta ,  \delta
\phi_1),
(d^+  _{\lambda} +r^*d^- _{\lambda}) \theta + \overline{\phi_2} \otimes  \delta
\phi_1,
\sD _{\lambda_1} \delta \phi_1 + \eta * \delta \phi_1 + \theta * \phi_2)
$$

 Assume that
$$
c_{dim \sM _{SW}} (Ind(d_{ \lambda} \oplus \sD _{ \lambda_1} )) \neq 0.
\tag3.7
$$
Assumption provides us with the nonzero point
$ (\theta ,  \delta \phi_1) \in
ker \delta_2 $
for some solution $a = \lambda_1 \oplus \lambda_2, \phi = 0 \oplus \phi _2 $
or in other words with a point in
$\widehat{ \sM _{B} ^{g, \omega}}(c+c_1(E), p_1(E))$ lying over
the corresponding Seiberg - Witten moduli space.

\proclaim{Lemma 3.1}
With the assumptions as above the differential (3.6) is epimorphic.
\endproclaim

\demo{Proof}
We shall remark first that if
$s_{\delta_2} (a, \phi_2, \langle \theta, \delta \phi_1 \rangle ) = 0$ for
$(a, \phi_2, \langle \theta, \delta \phi_1 \rangle ) \in
\Bbb P( ker \delta_1 ^* )$
then
$$
(\theta , \delta \phi_1) \ne 0 \Rightarrow \theta \ne d_{\lambda} \xi \and
\delta \phi_1 \ne 0.
$$
Indeed, assume $\theta = d_{\lambda} \xi$. Vanishing of the section implies in
particular
$$
d_{\lambda} ^* \theta + m^* _{\phi _2} (\delta \phi_1) = 0, ~~~
d_{\lambda} ^+ \theta + (\delta \phi_1 \otimes \phi_2)_0 = 0.
$$
This is equivalent to
$$
d_{\lambda} ^* d_{\lambda} \xi + (\delta \phi_1 \otimes \phi_2)_{tr} = 0
\tag3.8
$$
$$
\sigma \otimes [u, \xi] +  (\delta \phi_1 \otimes \phi_2)_0 = 0.
\tag3.9
$$
Taking scalar product of (3.8) with $\xi$ one has
$$
(\sigma \otimes [u, \xi], \xi) +  ((\delta \phi_1 \otimes \phi_2)_0, \xi) =
0+ ((\delta \phi_1 \otimes \phi_2)_0, \xi)=0
$$
and small calculation shows that
$$
 ((\delta \phi_1 \otimes \phi_2)_0, \xi)=0 \Rightarrow
 ((\delta \phi_1 \otimes \phi_2)_{tr}, \xi)=0
$$
Now taking scalar product of (3.9) with $\xi$ one has
$$
(d_{\lambda} ^* d_{\lambda} \xi, \xi) + ((\delta \phi_1 \otimes \phi_2)_{tr},
\xi) =
(d_{\lambda} ^* d_{\lambda} \xi, \xi) = 0 \Rightarrow \vert d_{\lambda} \xi
\vert = 0.
\Rightarrow \xi = 0.
$$
Substituting this to (3.8) and (3.9) we get
$(\delta \phi_1 \otimes \phi_2) = 0$
which, provided $\phi_2$ being a solution to an elliptic equations with
Laplace - type symbol vanishes totally if vanishes
on the open subset, means that $\delta \phi_1 = 0$. This is a contradiction to
the assumption
$(\theta , \delta \phi_1) \ne 0$ .

In a similar way let us assume $\delta \phi_1 = 0$. From the vanishing of the
section it follows
 in particular
$$
 \sD_{\lambda_1} ( \delta \phi_1)
  +(\theta) \cdot \phi_2 =0
$$
which means $(\theta) \cdot \phi_2 =0$ and therefore $ \theta = 0$ for the same
reason as above.

Take the nonzero point
$ (\theta ,  \delta \phi_1) \in
ker d_{ \lambda} \oplus \sD _{ \lambda_1}  $.
Then, as we showed, $ \theta \neq  d_{\lambda} \xi$ and by Thm 4.19 [FU]
restriction of the differential of our map
to infinitesimal variations of metrics $g=g+r$ covers
$  \Omega ^2 _+ (L_1 \otimes L_2 ^{-1})$ : $ r \to r^* d^- _{ \lambda} \theta$
is onto
$coker d^+ _{\lambda}$.
As well $   \delta \phi_1 \neq 0$ and infinitesimal variation of 1-forms
covers $\Gamma (W^- \otimes L_1)$:
$ \eta \to \eta * \delta \phi_1$ is onto $coker \sD _{\lambda_1} $.
Together with Lemma 1.1 this proves the statement.
\enddemo

Therefore one can assume that for generic metric and 1-form
second homology group vanishes and reducibles of the second type
are in general position.
Thus for a generic parametra one has zero locus of the section $s_{\delta_2}$
being
 smooth finite-dimensional submanifold of
the mentioned projective fibration with the fibre
$\Bbb P( \Omega ^1 (L_1 \otimes L_2 ^{-1})
  \oplus \Gamma (W^+ \otimes L_1) / \delta_1 (\Omega ^0 (L_1 \otimes L_2
^{-1})))$ over the
point of $\sM_{SW} (\beta)$. This submanifold will be referred to as a link of
$\sM_{SW} (\beta)$ in $\widehat{ \sM _B}$.

Following theorem shows that when some extra continuous parametra are added
this is true also for non reducible solutions.
As above dependence on the choice of metric and form $\omega$ gives gauge
invariant smooth map:
$$
\sT \times \sM etr \times \Omega^1 \times \sA_{\omega} \times
\Gamma (W^+ \otimes E) @>w \oplus v>>
\Omega^2_+ (pu(E)) \times \Gamma (W^- \otimes E)
\tag3.10
$$
defined as
$$ (t, g, \omega , a, \phi) \to (\sD_a \phi, F_a ^{+_g}
 (t-1)((- \phi \otimes \overline{\phi}))_{00}
\tag3.11
$$

If we prove that it is a submersion it will follow from Sard-Smale theorem that
for generic $g$ and $\omega$ the moduli space is smooth in the sense that the
second homology of corresponding deformation complex vanishes.

\proclaim{Theorem 3.1} The map $w$ is a submersion.

\endproclaim
\demo{Proof} Let $Dw, Dv$ denote a differentials of maps $w, v$ resp. It
follows
from [P-T, Prop.1.3.5] that if $D_{(a,\phi)}v_{\vert \delta a = \delta g
=\delta t = 0}$
isn't onto
then the connection $a$ is reducible $a = \lambda_1 \oplus \lambda_2$ and the
spinor has form $ \phi = 0 \oplus \phi_2$, that is the case of the previous
lemma.
 Now we shall use variations of metric and connection to prove that
$D_{(a,\phi)}w_{\vert \delta \phi = \xi = 0} = \Omega^2_+(pu(E))$.

Condition
$$
 \langle d_a^+(\delta a) + r^* F^-_a +
(\delta t ( \phi \otimes \overline{\phi}))_{00}
, \Phi \rangle = 0
\tag3.12
$$
for self dual non zero $\Phi \in \Omega^2_+(u(E)) $ is equivalent to
$$
d_a(\Phi) = 0 = d_a^*(\Phi)
\tag3.13
$$
and images of $F^-_a$ and $\Phi$ considered as map from
$\Lambda^2_-, \Lambda^2_+$ resp.  to $pu(E)$ are orthogonal as well as
those of $( \phi \otimes \overline{\phi})$ and $\Phi$
(cf. [F-U], Lemma 3.7). That is $\Phi$ does not vanish on an open set.
Thus maximal possible rank of the image of $F^-_a$ is 2. In that case $\Phi$
has rank one in generic point and therefore may be written as
$\Phi = \chi \otimes u$, where $\chi \in \Omega^2_+, u\in \Omega^0(u(E))$
and $\vert u \vert = 1$. Standard computation as in [F-U, Thm 3.4] shows that
 $d \chi = d_a u = 0$. Thus image of $F^-_a$ is perpendicular to $u$.
 The same reference gives $(F^-_a, u) = 0 \Rightarrow [F^-_a,u] \neq 0$.
So $ [F_a, u] = [F^+_a, u] \oplus [F^-_a, u] \neq 0$ and
$$
d_a u = [F_a, u] = [F^-_a, u] \oplus [-\phi \otimes \overline{\phi}, u] \neq 0
\tag3.14
$$
which gives a contradiction. This means that $F^-_a$
 and $(- \phi \otimes \overline{\phi})$
has rank at most one and it is exactly one at generic point
and their images are parallel.
This is possible if only $\phi $ considered as a map $W^+ \to E$ has rank one.
In this case by [PT, Prop. 1.3.5]
$a$ is reducible as well and we are done.

\enddemo

Results of this section can be formulated as follows:

\proclaim{Theorem 3.2}
For generic choice of parameter in
$\sT \times \sM etr \times \Omega^1 \times$
moduli space $\widehat{ \sM _{B} ^{t, g, \omega}}(c+c_1(E), p_1(E))$
is a smooth manifold of dimension given by (2.16) with the boundary
$$
\Cal {MP}  ^{g, \omega}(c+c_1(E), p_1(E)) \cup
\bigcup_{\beta \vert (f-\beta)^2 = p_1} \text{Link of}~ \sM_{SW} (\omega^+,
t,\beta)
$$
\endproclaim

\proclaim{Remark}
In the case $b_2 ^+ (X) =1$ one has chamber structure for both invariants.
Specifying chambers for Seiberg-Witten invariants related to some fixed
chamber for Donaldson polynomial one has to take in account shift of the
period by $\omega ^+$ which may a priory move it out of the positive cone
(i.e. period space of the Donaldson theory).
\endproclaim

\head
\S 4. Compactification.
\endhead

The Weitzenb\"ock formula
$$
\sD_a ^* \sD_a \phi = \nabla _a ^* \nabla _a \phi - (F_{a_0} ^+  +
F_{\nabla_{det}}) \phi +
\frac s 4 \phi
$$
provides universal estimates:
$$
\vert \phi \vert \le const_1
$$
$$
\vert F_a ^+ \vert ^2 \le const_2
\tag4.1
$$
in the same way it does for $U(1)$ case (compare (1.10)) although in this case
it depends not
only on the scale curvature of the manifold but also on the maximum of the
curvature $F_{\nabla_{det}}$:
$$
0 \le \Delta  \vert \phi \vert \le \frac {-s} 4  \vert \phi \vert ^2
+ \langle ((t-1)(\phi \otimes \overline {\phi}))_{00} (\phi), \phi \rangle
+  \langle (F_{\nabla_{det}} (\phi), \phi \rangle
=
$$
$$
= \frac {-s} 4  \vert \phi \vert ^2 - \frac 1 2 \vert \phi \vert ^4 +
+ \langle (t(\phi \otimes \overline {\phi}))_{00} (\phi), \phi \rangle
+  \langle (F_{\nabla_{det}} (\phi), \phi \rangle \le
$$
$$
\le (\frac {-s} 4 + \vert F_{\nabla_{det}} \vert) \vert \phi \vert ^2 -
 ( \frac 1 2 - \vert t \vert) \vert \phi \vert ^4 \le
(\frac {-s} 4 + \vert F_{\nabla_{det}} \vert) \vert \phi \vert ^2 -
 \frac 1 4 \vert \phi \vert ^4
$$
provided $ \vert t \vert \le \frac 1 4 $ (i.e. $t \in \sT _{\frac 1 4}$).
Therefore
$$
\vert \phi \vert \le -s + 4 \vert F_{\nabla_{det}} \vert.
\tag 4.1'
$$

 Weil formula for $c_2(E)$ gives
$$
 8 \pi ^2 c_2(E)\le \| F^- _a  \| _{L_2} \le 8 \pi ^2 c_2(E) + const_2 \cdot
vol_X.
\tag4.2
$$
This inequality provides following
\proclaim{Lemma 4.1}
Let $(a_i, \phi_i)$ be a sequence of solutions of our system and
$\epsilon $  - a real positive number.
Then there is a subsequence $i_j$ and a finite set $x_1, ..., x_p \in X$
such that
$$
\forall y \in X - \{x_i \} ~ \exists D_y \owns y \text{ such that} ~
\forall j ~ \| F_{a_{i_j}} \| _{L_2} ^2 \le \epsilon.
$$
\endproclaim

By this lemma and Theorem 2.3.7 of [DK] for each such $D_y$
there exist a gauge in which one has an estimate
$$
\| a_{\vert D _y} \| _{L_1 ^2} \le M
\| F_{a_{\vert D _y}} \| _{L ^2}
$$
for some constant $M$.

To get $C^{\infty} $ - convergence on the disk $D_y$ (or in fact on the smaller
disk
$ \hat{D}_y  \Subset D_y$) for any $y \in X - \{x_i \}$, and therefore on a
punctured manifold $X - \{x_i \}$, one uses approach of [DK], Ch 2, for
the proof of Uhlenbeck theorem.

The only place one needs antiselfduality
for the proof of the Theorem 2.3.8 of [DK] is getting from $L ^2 _1$ bound on
connection a $L_2 ^2$ bound on smaller domain.
This one has to be replaced by an estimate:
$$
\| a_{\vert \hat{D} _y} \| _{L_2 ^2} \le
const \cdot (\| a_{\vert D_y } \| _{L_1 ^2}
+ \| \phi \| _{L_1 ^2}) \le
const \cdot \| a_{\vert D_y } \| _{L_1 ^2} + const_2 \cdot vol _{D_y}
\tag4.3
$$
where $ \hat{D}_y  \Subset D_y$ which gives at the end an estimate
$$
\| a_{\vert \hat{D} _y} \| _{L_l ^2} \le M_{l,D_y} \| F_a \| _{L^2 (D_Y)}
+ const \cdot vol(D_y)
\tag4.3'
$$
while the rest of the scheme applied is unchanged.  The estimate (4.3')
provides
$C^{\infty} $ -convergence on any compact subset of $X - \{x_i \}$ together
with an
estimate of the norm of curvature of the limit $(a_{\infty}, \phi _{\infty})$:
$$
\| F_{a_{\infty}} \| _{L^2} \le 8 \pi ^2 c_2(E) + const_2 \cdot vol_X
$$

Removing of singularities is also a slight modification of that for
asd-connections.
Since harmonic spinor is smooth (provided we coupled the Dirac operator with a
smooth connection $a$)
 one has regularity of $L_1 ^2$ solutions of our system.
Cutting off the singularity $(\psi^2 a_0, \psi \phi)$ in a small disk $D_x$
centred in a singularity
with a smooth  cut-off function $\psi$ vanishing in a $r$-disk and equal to one
outside of
$2r$-disk gives a small
error to the first equation of the system:
$$
F_{\psi^2 a_0} +(\psi \phi \otimes \psi \phi)_{00} =
\psi ^2 (F_{a_0} + (\phi \otimes \phi)_{00}) + d(\psi ^2)\cdot a_0 +
( \psi ^4 - \psi ^2) [a_0, a_0]^+
$$
$$
\| F_{\psi^2 a_0} +(\psi \phi \otimes \psi \phi)_{00} \|_{L^2(D_x)} \le
\| d(\psi ^2) \| _{L^4(D_x)} \| a_0 \| _{L^4(D_x)}+ \| a_0 \| _{L^4(D_x)} \le
$$
$$
\le const \| F_{a_0} \| _{L^2(D_x)},
$$
with a similar estimate for $\| \sD _{\psi a} (\psi \phi) \| _{L^2}$.
Taking weak $L^2_1$ limits in a Coulomb gauge when $r$ tends to zero gives
a limit pair $(a', \phi')$ with smooth connection and spinor
provided by the above mentioned regularity of $L_1 ^2$ solutions.

Let us add ideal points to our moduli space. These are triples:
$$
((a, \phi), (x_1, ... , x_l)), ~~~ \text{where} ~~~ (a, \phi) \in \sM_B (p_1
+4l),
x_i \in X.
\tag4.4
$$
We shall assign a (curvature) density
$$
\vert F_a \vert ^2 +8 \pi \sum_{i=1} ^l \delta _{x_i}
\tag4.5
$$
to each of these and say that a sequence of points $(a_j, \phi_j)$
converges to the ideal point if curvature densities converges as measures
and there is a $C ^{\infty}$ convergence on compact subsets in
$X \backslash \{ x_j \} $. Let us omit for simplicity some indexes:

$ \sM _{B} ^{g, \omega}(c+c_1(E), p_1(E)) =  \sM _{B}( p_1)$.
\proclaim{Definition 4.1}
Let $\overline { \sM _{B}( p_1)} $ be a closure of the moduli space in the
union
$$
\sM_{B}(p_1) \cup \sM_{B}(p_1 +4) \times X \cup ...
\cup \sM_{B}(p_1+4k) \times s^kX \cup ...,
\tag4.6
$$
endowed with the convergence topology.
\endproclaim

As a result one has following
\proclaim{Theorem 4.1}
The space  $\overline {\sM_{B}} $ is compact
\endproclaim

Dimension counting shows that for generic parametra
only a subvariety of codimension
$2k$ of the stratum
$$
\sM_{B}(p_1+4k) \times s^kX
\tag4.7
$$
is in the closure of  $\sM_{B}$.

\proclaim{Corollary}
If $m_i \in \sM_{B} (p_1) \cap V_{\Sigma}$ is a sequence converging
to a point $(m_{\infty}, \{ x_1, ..., x_k \}) $ of the compactification, then
either $x_i \in \Sigma$ or $m_{\infty} \in \sM_{B} (p_1 + 4k) \cap V_{\Sigma}$.
\endproclaim

\head
\S 5. Localisation of the Donaldson polynomial.
\endhead

Now one can look at our moduli spaces as at bordism connecting
links of reducibles points. Reducibles of the first type,
moduli space $ \Cal {MP}$, is a boundary of $ \sM _{B}$ and its link
is a copy of $\Cal {MP}$ itself (or relevant compactification if
considering $\overline { \sM_{B}}$). The set of reducibles of the
second type is union of all Seiberg-Witten moduli spaces $\sM_{SW} (\beta)$
such, that  $ (- \beta + f)^2 = p_1$ or $ ( - \beta +f)^2 \ge p_1$
if one takes in account compactification $\overline { \sM_{B}}$.

We shall take a Poincare dual to certain cohomology classes
on our moduli space $\sM _{B}$ and on its compactification
$\overline { \sM_{B}}$.
If one takes Poincare dual on $\overline { \sM_{B}}$ to the
cohomology class of degree $dim \sM_B - 1$ one has a 1-dimensional singular
manifold $\sI$ with singularities diffeomorphic to the real
cone over finite number of points.
We shall assume that this manifold intersect lower strata
of the compactification only at reducible points of the
second type. Indeed dimension of the moduli space drops by 6:
$$
dim \sM_{B} (p_1) - dim \sM_{B} (p_1 + 4) = 6.
\tag5.1
$$
Therefore if one bubbling point ``kills'' not more then
two degree 2 classes (like $\mu_{\Sigma}$ or $t$) or one degree
4 class  ( $\nu$ -class) then $\sI$ intersects lower
strata of the compactification in reducible points only.
This works since one takes representatives ``localised''
over Riemannian surfaces (for $\mu_{\Sigma}$) or a point
(for $\nu$ or $t$) unless there is a trivial connection
in the compactification.
Since all $ \mu $-classes are lifted from the moduli of
$asd$-connections one can use standard dimension counting
to show that $\sI$ avoids lower strata in the compactification
of reducibles of the first type
$ \overline {\Cal {MP}}$
 and the number $\sharp \partial \sI \cap \Cal {MP}$
is the value of the relevant Donaldson polynomial (or Spin polynomial
more generally).
Thus manifold with boundary $ \sI$  is relating certain
intersection number on $\Cal {MP}$
with certain integrals over links of singularities of the
second type in $\overline { \sM_{B}}$. If one takes, say,
class $\prod _1 ^d \mu_{\Sigma} \times t^n$,
$2d = v.dim \sM_{asd} (p_1) = -2p_1 - 3(1+b_2^+)$
one has following equality:

$$
\gamma ^d _{p_1, f \text{\it mod2} - w_2(X)} ( \Sigma ) =
\sum _{ ( - \beta +f)^2 \ge p_1} \int _{ \text {Link of} \sM_{SW} (\beta)}
 \prod_1 ^d \mu_{\Sigma} u^n.
\tag5.2
$$
Of course the summand in the r.h.s is nontrivial only if
the moduli space $\sM_{SW} (\beta)$ is nonempty.

Now to describe a value of the Donaldson polynomial we have to describe for any
Seiberg-Witten class $\beta$ the link of
$$
\sM_{SW}(\beta) \subset Sing_2 \subset \sM_{B}(f, p_1)
\tag5.3
$$
( We omit the upper indexes of the moduli space because it doesn't depend on
the continuous
parametra) and the integral of the type (5.2) over this link. The point is that
this
integral is the standard polynomial of the linear forms
$$
\langle \beta,  \rangle, \langle f,  \rangle
$$
and the quadratic intersection form $q_X$. Let call it the local polynomial and
denote it as $loc \gamma
$.

Now in the localisation formula (5.2) the left side doesn't depend on $f$ and
we can consider one as the relations between the local polynomials. Using this
relations we will compute the local polynomials precisely in the next paper.

For this computation we will construct pure geometrical finite space
$MH^1(\beta, r)$, where $4r = (-\beta + f)^2 - p_1$ is number of points on $X$,
as the space of "gluing parameters", (which doesn't depend on 2-cohomology
class $f$) and the virtual vect
or bundle on  $MH^1(\beta, r)$ such that the intersections of $\mu$-classes
with "top" Chern class of this bundle give the local polynomials. These
constructions are pure geometrical and here we can say "goodbye" to the gauge
theory,
connections , moduli spaces and all standard stuff of the Donaldson's Theory.

\enddocument